\documentclass[fleqn,11pt]{article}
\usepackage{latexsym}
\usepackage{graphicx}
\usepackage{amsmath}
\usepackage{amssymb}
\usepackage{amsfonts}
\usepackage{epsfig}
\usepackage{color}
\newcommand{\qed}{\hfill$\Box$\par\smallskip\noindent}

\newtheorem{definition}{Definition}[section]

\newtheorem{theorem}{Theorem}[section]
\newtheorem{lemma}[theorem]{Lemma}
\newtheorem{proposition}[definition]{Proposition}

\newcommand{\comment}[1]{{}}

\begin{document}



\title{The Largest-Z-ratio-First algorithm is 0.8531-approximate for \\
scheduling unreliable jobs on $m$ parallel machines
}


\author{
Alessandro Agnetis \thanks{Dipartimento di Ingegneria
dell'Informazione, Universit\`a di Siena, e-mail
agnetis@dii.unisi.it} \and Thomas Lidbetter \thanks{ Management Science and Information Systems Department, Rutgers Business School, Rutgers University, e-mail
tlidbetter@business.rutgers.edu} }

\date{}
\maketitle
\medskip \noindent

\begin{abstract} \noindent
In this paper we analyze the worst-case performance of a greedy algorithm called Largest-Z-ratio-First for the problem of scheduling unreliable jobs on $m$ parallel machines. Each job is characterized by a success probability and a reward earned in the case of success. In the case of failure, the jobs subsequently sequenced on that machine cannot be performed. The objective is to maximize the expected reward. We show the algorithm provides an approximation ratio of $\simeq 0.853196$, and that the bound is~tight.
\end{abstract}


\section{Introduction}\label{sec:intro}
The following problem, called \emph{Unreliable Job scheduling Problem} (UJP), was introduced in \cite{ADPS:jos}. A set of jobs $J=\{J_1, \ldots, J_n\}$ must be assigned to $m$ parallel, identical machines, $M_1, \ldots, M_m$. Each job must be assigned to a single machine and a machine can process one job at a time. Jobs are \emph{unreliable}, i.e., while a job is being
processed by a machine, a \emph{failure} can occur, which implies losing all the work which was scheduled but not yet executed by the
machine. Each job $J_i$ is characterized by a certain \emph{success probability} $\pi_i$ (independent from other jobs) and a
\emph{reward} $r_i$, which is gained if the job is successfully completed. We assume that the values $\pi_i$ are rational numbers strictly between 0 and 1. (A job with probability 0 will always be appended to the end of a machine schedule, without affecting the value of the objective function, while a job with probability 1 can always be scheduled first on a machine, only adding a constant term to the objective but without affecting the rest of the schedule.) The problem is to find an assignment of jobs to the $m$ machines and a sequence on each machine that maximizes the total expected reward. We may also frame UJP as a search and rescue problem, as in~\cite{Lidbetter}. 

In this note we consider the Largest-Z-ratio-First algorithm for UJP, and, for any instance $I$, give a bound on the ratio $\lambda(I)$ between the value produced by the Largest-Z-ratio-First algorithm and the optimal value. Limited to the case $m=2$, it was shown in \cite{ADP2014} that $\lambda(I)\geq (2+\sqrt{2})/4\simeq 0.8535...$. Here we extend the result to any value of $m$, showing that $\lambda(I)\geq 0.85319...$.

In Section \ref{sec:review} we review some basic notions concerning UJP, while in Section \ref{sec:proof} we provide the main result.

\section{Unreliable Job scheduling Problem}\label{sec:review}

Here we briefly review the main concepts and notation concerning UJP. Let $S^h$ be a
sequence of $K$ jobs assigned to machine $M_h$, and let $S^h(k)$ be the job in $k$-th position in $S^h$. A feasible solution
$S=\{S^1,S^2,\dots,S^m\}$ for UJP is an assignment and sequencing of the $n$ jobs on the $m$ machines. If $K$ jobs are assigned to $M_h$, the \emph{expected reward} of sequence $S^h$ is given by
\begin{equation}\label{funzob}
ER[S^h]=\pi_{S^h(1)}r_{S^h(1)}+\pi_{S^h(1)}\pi_{S^h(2)}r_{S^h(2)}
+\ldots+\pi_{S^h(1)}\ldots
\pi_{S^h(K-1)}\pi_{S^h(K)}r_{S^h(K)}.
\end{equation}
and the total expected reward is therefore
\[
ER[S]=ER[S^1]+ER[S^2]+ \ldots+ER[S^m].
\]

UJP consists of finding a solution $S_{OPT}=\{S_{OPT}^1,S_{OPT}^2,\ldots,S_{OPT}^m\}$ that maximizes the total expect reward.
A key role in our analysis is played by the following quantity associated with each job $j$, called the \emph{Z-ratio}:
\begin{equation}\label{Zrat}
Z_j=\frac{\pi_jr_j}{1-\pi_j}.
\end{equation}
When $m=1$, the optimal solution is achieved by sequencing the jobs in non-increasing order of $Z_j$ \cite{mitten,ADPS:jos}. Hence, UJP indeed consists of deciding how to partition the $n$ jobs among the $m$ machines, since on each machine the sequencing is then dictated by the priority rule (\ref{Zrat}). Since  \textsc{Product Partition} can be polynomially reduced to UJP with $m=2$ \cite{ADPS:jos}, and since \textsc{Product Partition} is strongly NP-hard \cite{ngetal}, so is UJP, even for $m=2$.

UJP bears various similarities with the classical  problem of minimizing total weighted completion time on $m$ parallel machines, i.e., $Pm||\sum_j w_jC_j$, though unlike this problem, UJP is not concerned with the processing times of jobs. The single-machine problem $1||\sum_j w_jC_j$ is solved by the well-known Smith's rule, i.e., sequencing the jobs in non-increasing order of the ratio $\rho_j=w_j/p_j$. For any $m\geq 2$, $Pm||\sum_j w_jC_j$ is NP-hard. 

The Largest-Ratio-First algorithm for $Pm||\sum_j w_jC_j$ is the following: order the jobs by non-increasing ratios $\rho_j$ and assign them in this order to the $m$ machines, allocating the next job in the list to the machine that frees up first.  A schedule obtained in this way is called a Largest-Ratio-First (LRF) schedule. Kawaguchi and Kyan \cite{kwkyan} showed that the worst-case error of any LRF schedule is $(1+\sqrt{2})/2$. A simpler proof of this result has been provided by Schwiegelshohn \cite{schw2011}.

In this paper we analyze the performance of an analogous ratio-based algorithm for UJP. In the following, while assigning the jobs to machines, we call the \emph{cumulative probability of a machine} the product of the success probabilities of the jobs already scheduled on that machine. When a job $j$ is assigned to a machine, we use the notation $P_j$ to indicate the product of the success probabilities of all jobs scheduled on the machine up to job $j$ (included), and we refer to $P_j$ as the \emph{cumulative probability of job $j$}.

The \emph{Largest-Z-ratio-First algorithm} for UJP works as follows. Order the jobs by non-increasing $Z_j$ and assign them in this order to the $m$ machines, allocating the next job in the list to a machine currently having maximum cumulative probability (ties are broken arbitrarily). A schedule obtained in this way is also called a \emph{Largest-Z-ratio-First} (LZF) schedule. In this paper we investigate the worst-case performance of any LZF schedule.

In establishing our result, we follow a similar line of reasoning to the one in \cite{schw2011} for $Pm||\sum_j w_jC_j$. While our Lemmas \ref{prop1} and \ref{prop2} are an adaptation of Corollaries 1 and 3 in \cite{schw2011}, Lemma \ref{prop3} is novel and exploits features that are specific to UJP.

\section{An approximation bound}\label{sec:proof}

The bound on the performance of an LZF schedule is proved by subsequently reducing the set of instances which need to be considered in order to detect the worst-case instance. This is done through three lemmas. In Lemma \ref{prop1} we show that we can restrict to instances in which all jobs have $Z_j=1$. In Lemma \ref{prop2} we prove that it is sufficient to consider instances containing at most $m-1$ jobs having a very large reward and an arbitrary number of jobs having small reward. Lemma \ref{prop3} shows that, furthermore, the worst-case situation occurs when the success probabilities of the jobs having very large reward are arbitrarily close to 0. Thereafter, the main result can be established.

As in \cite{schw2011}, we extend the usual definition of an instance $I$ of the problem to include an arbitrary LZF order for all the jobs. This order produces the \emph{primary LZF} schedule $S_{LZF}(I)$. In this way, an instance $I$ has a unique primary LZF schedule, even if all jobs have the same $Z_j$. For an instance $I$, we let
\[
\lambda(I)=\frac{ER[S_{LZF}(I)]}{ER[S_{OPT}(I)]}. \]

The following lemma is an adaption of Corollary 1 in \cite{schw2011}.
\begin{lemma}\label{prop1}
For every instance $I$ of UJP, there is an instance $I'$ with $\lambda(I')\leq \lambda(I)$ and $Z_j=1$ for all jobs $j\in I'$.
\end{lemma}
Proof. Let $\zeta_1 > \zeta_2 >\dots>\zeta_d$ be the $d$ different $Z_j$ values of jobs in $I$, and let $\zeta_{d+1}=0$. We can write $\zeta_i$ as $\zeta_i=\sum_{k=i}^d(\zeta_k-\zeta_{k+1})$. Letting $i(j)$ denote the index of the Z-ratio of job $j$, we have
\[
r_j=\zeta_{i(j)}\frac{1-\pi_j}{\pi_j}=\sum_{k=i(j)}^d(\zeta_k-\zeta_{k+1})\frac{1-\pi_j}{\pi_j}.
\]
Recalling the definition of $P_j$, the expected reward of an arbitrary schedule $S$ for instance $I$ can be written as
\[
ER[S]=\sum_{j\in I} r_jP_j=\sum_{j\in I} \left(\sum_{k=i(j)}^d(\zeta_k-\zeta_{k+1})\frac{1-\pi_j}{\pi_j}\right)P_j=
\]
\begin{align}
=\sum_{k=1}^d\left((\zeta_k-\zeta_{k+1})\sum_{j:i(j)\leq k}\frac{1-\pi_j}{\pi_j} P_j\right). \label{eq:ER}
\end{align}

Next, we define a sequence of instances of UJP, $I_k=\{j\in I:i(j)\leq k\}$, $k=1,\dots,d$. For these instances, we set the reward of job $j$ in $I_k$ as $(1-\pi_j)/\pi_j$ (hence, $Z_j=1$). It follows that {\em any} ordering of the jobs in $I_k$ is an LZF order, so we can select an LZF order for $I_k$ that is consistent with our LZF order for $I$, ensuring that for any job $j\in I_k$, the values of $P_j$ in $S_{LZF}(I_k)$ and in $S_{LZF}(I)$ are identical. Letting $P_j^*$ denote the value of $P_j$ in an optimal schedule $S_{OPT}(I)$, and observing that $ER[S_{LZF}(I_k)] = \sum_{j \in I_k}((1-\pi_j)/\pi_j) P_j$, we can apply~(\ref{eq:ER}) to both $S_{LZF}(I)$ and $S_{OPT}(I)$ to obtain
\[
\lambda(I)=\frac{\sum_{k=1}^d(\zeta_k-\zeta_{k+1})ER[S_{LZF}(I_k)]}{\sum_{k=1}^d(\zeta_k-\zeta_{k+1})\sum_{j\in I_k}\frac{1-\pi_j}{\pi_j}P_j^*}.
\]
By the optimality of $S_{OPT}(I_k)$, we must have $ER[S_{OPT}(I_k)]\geq \sum_{j\in I_k}((1-\pi_j)/\pi_j)P_j^*$, so that
\[
\lambda(I)\geq \frac{\sum_{k=1}^d(\zeta_k-\zeta_{k+1})ER[S_{LZF}(I_k)]}{\sum_{k=1}^d(\zeta_k-\zeta_{k+1})ER[S_{OPT}(I_k)]}\geq \min_{1\leq k \leq d}\lambda(I_k).
\]
Hence, $\lambda(I)$ is at least as large as the value it attains in an instance in which all jobs have Z-ratio equal to 1.\qed

Notice that if $Z_j=1$ for all jobs, the expected reward of the jobs scheduled on a certain machine $M_h$, from (\ref{funzob}), is given by

\[
ER[S^h]=\pi_{S^h(1)}\left(\frac{1-\pi_{S^h(1)}}{\pi_{S^h(1)}}\right)+\pi_{S^h(1)}\pi_{S^h(2)}\left(\frac{1-\pi_{S^h(2)}}{\pi_{S^h(2)}}\right)
+\ldots+
\]
\begin{equation}\label{funzobz}
+\pi_{S^h(1)}\ldots
\pi_{S^h(K-1)}\pi_{S^h(K)}\left(\frac{1-\pi_{S^h(K)}}{\pi_{S^h(K)}}\right)= 1-\prod_{i=1}^K \pi_{S^h(i)}
\end{equation}

and hence, given a schedule $S$, if $P_h(S)=\prod_{i=1}^K \pi_{S^h(i)}$ is the cumulative probability of machine $M_h$ in schedule $S$, the expected
reward $ER[S]$ is given by
\begin{equation}\label{totER}
ER[S] = m-\sum_{h=1}^m P_h(S).
\end{equation}

In view of Lemma \ref{prop1}, from now on we only consider instances of UJP in which all jobs have Z-ratio equal to 1.

Given an instance $I$ and the corresponding primary schedule $S_{LZF}$, let $P_{\max}(I)=\max_{1\leq h\leq m}\{P_h(S_{LZF}(I))\}$. We can now establish the UJP counterpart of Corollary 3 in \cite{schw2011}.

\begin{lemma}\label{prop2}
For every instance $I$ of UJP for which all jobs have Z-ratio equal to 1, there is an instance $I'$ such that $\lambda(I')\leq \lambda(I)$ and every job has an arbitrarily high success probability if its cumulative probability in $S_{LZF}(I')$ is at least $P_{\max}(I')$.
\end{lemma}
Proof. Consider an arbitrary  instance $I$ and the corresponding LZF schedule $S_{LZF}(I)$. Now consider an instance $J$ obtained by replacing any job $j$ with two jobs $j1$ and $j2$ such that $\pi_{j1}\pi_{j2}=\pi_j$, and consider the schedule $S(J)$ obtained from $S_{LZF}(I)$ replacing $j$ with $j1$ and $j2$ consecutively scheduled in this order on the same machine. Note that $ER[S_{LZF}(I)]=ER[S(J)]$, by~(\ref{totER}). Call $\bar P$ the cumulative probability of the jobs \emph{preceding} $j$ on the same machine in $S_{LZF}(I)$. Due to the mechanism of the LZF algorithm and since $\pi_i<1$ for all jobs, $\bar P> P_{\max}(I)$. We choose $\pi_{j1}$ large enough so that $S(J)$ is still an LZF schedule and therefore $\pi_{j1}\bar P\geq P_{\max}(I)$. (This is not the case if $\pi_{j1}\bar P< P_{\max}(I)$, as an LZF schedule would have assigned $j2$ on the machine that in $S_{LZF}(I)$ has cumulative probability $P_{\max}(I)$.) Also, note that $ER[S_{OPT}(J)]\geq ER[S_{OPT}(I)]$. Hence, as long as $S(J)$ is a LZF schedule,
\[
\lambda(I)=\frac{ER[S_{LZF}(I)]}{ER[S_{OPT}(I)]}\geq \frac{ER[S(J)]}{ER[S_{OPT}(J)]}=\lambda(J).
\]
We can repeat this job splitting until all jobs $j$ such that $P_j\geq P_{\max}(I)$ have a success probability that is arbitrarily close to $1$. Note all the jobs such that $P_j< P_{\max}(I)$ are the last scheduled jobs on the various machines, and hence there can be at most $m-1$ of them. \qed

The consequence of Lemma \ref{prop2} is that we can restrict to instances satisfying the following. Each machine processes a large number of jobs with an arbitrarily high success probability, until its cumulative probability reaches $P_{\max}(S_{LZF}(I))$. We call these jobs \emph{low-value jobs}, since, if the success probability is high, then the reward is low. After these jobs, at most $m-1$ machines process one more job. We call these jobs \emph{second-stage jobs}. So, in summary, from now on we only consider instances which contain several low-value jobs followed by at most $m-1$ second-stage jobs. 

We now take the last fundamental step, which consists of showing that the most unfavourable situation occurs when the success probabilities of all second-stage jobs are arbitrarily close to $1$ or arbitrarily close to 0. The proof of this lemma uses arguments that are specific to UJP and are not derived from \cite{schw2011}. 

\begin{lemma}\label{prop3}
 For every $\varepsilon > 0$ and every instance $I$ of UJP for which all jobs have Z-ratio equal to 1, there is an instance $I'$ such that $\lambda(I')\leq \lambda(I)$ and every job has an arbitrarily high success probability if its cumulative probability in $S_{LZF}(I')$ is at least $P_{\max}(I')$. All other jobs have success probability greater than $1-\varepsilon$ or less than $\varepsilon$.
\end{lemma}
 Proof. Given an instance $I$, let $p = P_{\max}(I)$ and denote by $t\in(0,m)$ the sum of all the success probabilities of the second-stage jobs in $I$ that lie in the range $[\varepsilon,1-\varepsilon]$. We will prove that we can construct an instance $I'$ with the required properties in the case that $t$ is an integer. In the case that $t$ is not an integer, let $q$ be an integer such that $tq$ is an integer. Given $I$, consider an instance $I^q$ in which each job of $I$ is replaced with $q$ consecutive copies, and there are $qm$ machines. In $I^q$, we adopt an LZF order producing $k$ identical copies of each machine schedule in $S_{LZF}(I)$, so that $ER[S_{LZF}(I^q)]=q ER[S_{LZF}(I)]$. It is also clear that $ER[S_{OPT}(I^q)] \ge q ER[S_{OPT}(I)]$, since any schedule for $I$ gives rise to a schedule for $I^q$ whose expected reward is $q$ times as large. Hence, $\lambda(I^q) \le \lambda(I)$, and it is sufficient to find a schedule $I'$ such that $\lambda(I') \le \lambda(I^q)$. But the sum of all the success probabilities of the second-stage jobs in $I^q$ that lie within $[\varepsilon,1-\varepsilon]$ is $tq$, hence an integer, so we may as well assume $t$ is an integer (otherwise, the whole argument is applied to $I^q$).

Recall that in $S_{LZF}(I)$ every machine processes a large number of low-value jobs with cumulative probability $p=P_{\max}(I)$, possibly followed by a second-stage job. Let $k=k(I)$ be the number of second-stage jobs in $I$ with success probability in the range $[\varepsilon,1-\varepsilon]$. For every fixed $m$, we will prove by strong induction on $k$ that there is an instance $I'$
such that $\lambda(I')\leq\lambda(I)$ and all the second-stage jobs of $\lambda(I')$ have success probability less than $\varepsilon$ or more than $1-\varepsilon$, which will prove the thesis.

This is evidently true
for $k = 0$, because in that case $I$ is already of the form $I'$. Note that we cannot have $k=1$, since in this case a single second-stage job would have success probability $\tilde\pi \in [\varepsilon, 1-\varepsilon]$ and $t=\tilde{\pi}$ would not be an integer.

Although not strictly necessary, we consider separately the case $k = 2$, since this is a
good introduction to the structure of the general induction argument. In this case, in $I$ there
must be two second-stage jobs $i$ and $j$ with success probabilities $\pi_i,\pi_j \in[\varepsilon,1-\varepsilon]$ such that $\pi_i+\pi_j=1$ (again, since $t$ is integer).
We define a new instance $I'$ by replacing jobs $i$ and $j$ with jobs $i'$ and $j'$ having success probabilities $\pi_i' = \delta$ and $\pi_j' = 1-\delta$ respectively, where $\delta\in(0,\varepsilon)$ will be specified later. This does not affect the expected reward because $\pi_i+\pi_j = \pi_{i'}+\pi_{j'} = 1$. Hence, the expected rewards of
$S_{LZF}(I)$ and $S_{LZF}(I')$ are the same, and $k(I')= 0$. We must show that there is a schedule
for $I'$ whose expected reward is at least $ER[S_{OPT}(I)]$, so that $\lambda(I')\leq\lambda(I)$.

Consider $S_{OPT}(I)$, and define a schedule $S'(I')$ for $I'$ which is the same as $S_{OPT}(I)$, but replacing $i$ and $j$ with $i'$ and $j'$, respectively. There are two possibilities. Either, in $S_{OPT}(I)$, jobs $i$ and $j$ are processed on the same machine or on different machines. If they are processed on the same machine, then the expected reward of $S'(I')$ must be greater than that of $S_{OPT}(I)$. This is because the
contribution of that machine to the expected reward of $S'(I')$ is $1-\hat p\delta(1-\delta)$, where $\hat p$ is the cumulative probability of all the other jobs processed by the machine; this is greater than the contribution $1-\hat p\pi_i\pi_j $ of the same machine to the expect reward of $S_{OPT}(I)$.
Now suppose that $i$ and $j$ are processed on different machines in $S_{OPT}(I)$, and let the cumulative probabilities of all the other jobs scheduled by $S_{OPT}(I)$ on these two machines be $p_1$ and $p_2$, respectively. First assume that $p_1\geq p_2$. Then the contribution of these two machines to the expected
reward of $S'(I')$ is $2-p_1\delta -p_2(1-\delta)$ and to $S_{OPT}(I)$ is $2-p_1\pi_i-p_2\pi_j$. Recalling that $\pi_i+\pi_j=1$, the difference between these two contributions is 
\[
p_1(\pi_i-\delta)+p_2(\pi_j-1+\delta)=(\pi_i-\delta)(p_1-p_2)\geq 0,
\]
as long as $\delta$ is chosen to be at most $\pi_i$.

If instead, $p_1 \leq p_2$, we define a different schedule for $I'$ which is the same as $S_{OPT}(I)$, but replacing $i$ with $j'$ and $j$ with $i'$. A similar argument to the one above shows that the expected reward of this schedule is at least $ER[S_{OPT}(I)]$ for $\delta \le \pi_j$.

Now consider any $k\geq 3$. Since we are using strong induction, we assume that the induction hypothesis is true for \emph{all} smaller values of $k$. Let $i$ and $j$ be any two second-stage jobs in $I$ with success probabilities $\pi_i,\pi_j \in [\varepsilon,1-\varepsilon]$. We consider two cases.

\emph{Case 1: $\pi_i+\pi_j\leq 1$}. This is similar to the case $k=2$. We define a new instance $I'$ by replacing jobs $i$ and $j$ with $i'$ and $j'$ that have success probabilities $\delta$ and $\pi_i+\pi_j-\delta$ respectively, where $\delta\in(0,\varepsilon)$ will be specified later. This does not affect the expected reward, so the
expected reward of $S_{LZF}(I)$ and $S_{LZF}(I')$ are the same. Since, by the induction hypothesis, we also have that $k(I') < k(I)$, it is sufficient to show that there is a schedule for $I'$ whose
expected reward is at least that of $S_{OPT}(I)$.
As before, we define a schedule $S'(I')$ for $I'$ which is obtained from $S_{OPT}(I)$ replacing
$i$ and $j$ with $i'$ and $j'$, respectively. Again, there are two possibilities: either $i$ and $j$ are processed on the same machine in $S_{OPT}(I)$ or not. In the former case, $ER[S'(I')]\geq ER[S_{OPT}(I)]$, by a similar argument as for $k=2$. In the latter case, let the cumulative probabilities
of all the other jobs scheduled by $S_{OPT}(I)$ on these two machines be $p_1$ and $p_2$, respectively. We assume that $p_1\geq p_2$. (If not, then as before, we define a different schedule for $I'$ in which $i$ is replaced with $j'$ and $j$ is replaced with $i'$.) Then the contribution of these two machines
to $ER[S'(I')]$ is $2-p_1 \delta - p_2(\pi_i+\pi_j-\delta)$ and to $ER[S_{OPT}(I)]$ is $2-p_1\pi_i-p_2\pi_j$. The difference between these contributions is 
\[
-p_1 \delta -p_2(\pi_i+\pi_j-\delta)+p_1\pi_i+p_2\pi_j=(\pi_i-\delta)(p_1-p_2)\geq 0,
\]
as long as $\delta$ is chosen to be at most $\pi_i$.

\emph{Case 2: $\pi_i+\pi_j> 1$}.
In this case, we define a new instance $I''$ by replacing jobs $i$ and $j$ in $I$ with jobs $i''$ and $j''$ that have success probabilities $\pi_i+\pi_j-1+\delta$ and $1-\delta$ respectively, where $\delta \in (0,\varepsilon)$ will be specified later. Again, from (\ref{totER}), $ER[S_{LZF}(I)]=ER[S_{LZF}(I'')]$, and $k(I'') < k(I)$, so by the induction hypothesis, it is sufficient to show that there is a schedule for $I''$ whose expected reward is at least that of $S_{OPT}(I)$. We define a schedule $S''(I'')$ for $I''$ which is obtained from $S_{OPT}(I)$ by replacing $i$ and $j$ with jobs $i''$ and $j''$ . Again there are two subcases. The first is when $i$ and $j$ are processed on the same machine in $S_{OPT}(I)$. In this case, let $\hat p$ be the cumulative success probability of all the other jobs processed on this machine. Then the contribution of this machine to the expected reward of $S''(I'')$ is $1-\hat p(\pi_i + \pi_j-1+\delta)(1-\delta)$ and to $S_{OPT}(I)$ is $1-\hat p\pi_i\pi_j$. The difference between the two contributions is therefore
\[
\hat{p}(\pi_i \pi_j-(\pi_i + \pi_j-1+\delta)(1-\delta))=\hat p((1-\pi_i)(1-\pi_j)-\delta(\pi_i+\pi_j+2-\delta)).
\]
Thus, if $\delta$ is chosen to be small enough, then the expression displayed above is positive.

The second subcase is when $i$ and $j$ are processed on different machines in $S_{OPT}(I)$. Again, let the cumulative
probabilities of all the other jobs scheduled by $S_{OPT}(I)$ on these two machines be $p_1$ and $p_2$ respectively, with $p_1\geq p_2$. Then the difference between the contributions of these two machines to $ER[S'(I')]$ and $ER[S_{OPT}(I)]$ respectively is
\[
(2-p_1(\pi_i+\pi_j-1+\delta)-p_2(1-\delta))-(2-p_1\pi_i-p_2\pi_j)= (1-\pi_j-\delta)(p_1-p_2),
\]
which is positive since $1-\pi_i-\delta>1-\pi_i-\varepsilon \ge 0$, and the proof is complete.\qed 

We are now in the position of establishing the main result of this paper.
\begin{theorem} \label{thm:main}
For any instance $I$ of UJP, $\lambda(I)\geq 0.853196...$
\end{theorem}
 Proof. By Lemmas \ref{prop1}, \ref{prop2} and \ref{prop3}, for any instance $I$ and any $\varepsilon >0$, we can find an instance $I'$ with $\lambda(I')\le \lambda(I)$ such that in $S_{LZF}(I)$, each machine processes some low-value jobs first, until the cumulative success probability
reaches some $p$ on each machine, then some machines process a second-stage job with success probability greater than $1-\varepsilon$ or less than $\varepsilon$. If machine $j$ processes a second-stage job, let $\pi_j$ be the success probability of that job; if it processes no second-stage job, let $\pi_j=1$. We may assume, by reordering, that $\pi_1,\ldots,\pi_t > 1- \varepsilon$ for some $t$ and $\pi_{t+1},\ldots,\pi_{m} < \varepsilon$. By choosing $\varepsilon$ to be small enough, we can ensure that in an optimal schedule for such an instance $I'$, for $j>t$, machine $j$ processes a single job with success probability $\pi_j$ and for $j \le t$, machine $j$ processes a set of jobs with cumulative probability $\pi_j p^{m/t}$. Recalling that $ER[S_{LZF}(I')]=m-p\sum_{j=1}^m \pi_j$, the ratio $\lambda(I')$ is given by
\[
\lambda(I')=\frac{m-p\sum_{j=1}^m \pi_j}{m-\sum_{j=1}^t \pi_j  p^{m/t}-\sum_{j=t+1}^m \pi_j}.
\]
Taking the limit as $\varepsilon \rightarrow 0$, we have $\pi_j\rightarrow 1$ for $j=1,\dots,t$ and $\pi_j\rightarrow 0$ for $j=t+1,\dots,m$, so that we obtain
\[
\lambda(I)\ge \frac{m-tp}{m-tp^{m/t}}=\frac{m/t-p}{m/t-p^{m/t}}=:f(m/t,p).
\]
The minimum of the function $f(m/t,p)$ (for $p\in(0,1)$ and $m/t \ge 1$) can be found by numerical methods and it has value 0.853195...{\footnote{This value and those in Table \ref{tab1} were confirmed by a number of nonlinear solvers, including Matlab with the GlobalSearch option.}.\qed}

The minimum of $f(m/t,p)$ is attained when $m/t\simeq 2.1231...$ and $p\simeq 0.58919...$.

The approximate values of $\min_{t,p} f(m/t,p)$ are shown in Table~\ref{table:min-f} for the first few values of $m$. In this table we also denote by $\bar t$ the value of $t$ for which the bound is attained.  We observe that for $m=2$ and $t=1$, simple calculus shows that the minimum of $f(2,p)$ is attained for $p=2-\sqrt{2}$ and its value is $(2+\sqrt{2})/4\simeq 0.85355...$, retrieving the result in \cite{ADP2014} for $m=2$.

\begin{table}[htb!]
\centering
\caption{Approximate minimum values of $f(m/t,p)$ and $\bar t$.}
\label{table:min-f}
\begin{tabular}{lll}
\\ \hline
$m$ & $\min_{t,p}f(m/t,p)$ & $\bar t$ \\
\hline
2 & $0.85355 $ &1\\
3 & $0.86179 $ &1\\
4 & $0.85355 $ &2\\
5 & $0.85541 $ &2\\
6 & $0.85355 $ &3\\
\hline
\end{tabular}\label{tab1}
\end{table}

(While the first few entries of the table hint at the possibility that the minimum of $f(m/t,p)$ is equal to $0.85355...$ for all even values of $m$, this is in fact not the case. For instance, when $m=38$, the minimum is $0.853199...$.) 

\begin{proposition}
 For each fixed $m$, an arbitrarily tight example can be built. 
 \end{proposition}
 
 Proof. Given $m$, let $\bar p$  and $\bar t$ be the values of $p$ and $t$ for which the minimum of  $f(m/t,p)$ is attained. Consider an instance in which $Z_j=1$ for all jobs. There are $m\bar t$ jobs with success probability $\pi$, where $\pi$ is a rational number arbitrarily close to $\bar p^{1/\bar t}$ (for the sake of simplicity, in what follows we consider $\pi=\bar p^{1/\bar t}$), and $(m-\bar t)$ jobs with success probability $\varepsilon$. (As $Z_j=1$, the rewards of the two job types are $(1-\pi)/\pi$ and $(1-\varepsilon)/\varepsilon$ respectively.)  Consider the schedule $\sigma_H$, in which $\bar t$ jobs with success probability $\pi$ are assigned to each machine, and $(m-\bar t)$ machines receive an additional job with success probability $\varepsilon$. Note that $\sigma_H$ is an LZF schedule, attained when first all jobs with probability $\pi$ are allocated, followed by all jobs with probability $\varepsilon$. Let $z_H$ be the expected reward of this schedule. In $\sigma_H$, the cumulative probability of $(m-\bar t)$ machines is $\varepsilon\bar p$, while that of the other $\bar t$ machines is $\bar p$. Hence, recalling \eqref{totER},$z_H=m-(m-\bar t)\varepsilon\bar p-\bar t\bar p$. 
 
 Now consider the schedule $\sigma^*$ in which $m$ jobs with success probability $\pi$ are scheduled on each of $\bar t$ of the machines, while the other $(m-\bar t)$ machines only receive a job with success probability $\varepsilon$. (Incidentally, we note that also $\sigma^*$ is an LZF schedule, obtained by scheduling first all jobs with probability $\varepsilon$ and then all the others.) So, in $\sigma^*$ the cumulative probability on each of the first $\bar t$ machines is $\bar p^{m/\bar t}$, while for the other $(m-\bar t)$ machines it is $\varepsilon$, thus yielding $z^*=m-(m-\bar t)\varepsilon-\bar t\bar p^{m/\bar t}$. We have:
\[
\frac{z_H}{z^*}=\frac{m-(m-\bar t)\varepsilon\bar p-\bar t\bar p}{m-(m-\bar t)\varepsilon-\bar t\bar p^{m/\bar t}},
\]
and letting $\varepsilon
\rightarrow 0$, we have
\[
\frac{z_H}{z^*}\rightarrow\frac{m-\bar t\bar p}{m-\bar t\bar p^{m/\bar t}}.
\]

As an example, consider the case $m=5$. In this case the value of the bound is 0.855411..., attained for $\bar t=2$ $\bar p=0.6024533...$. We define an instance with $m\bar t=10$ jobs having probability $\pi=\bar p^{1/\bar t}=0.6024533^{1/2}=0.776179$, and $(m-\bar t)=3$ jobs having probability $\varepsilon$. It holds
\[
5-\bar t\bar p=5-2^*0.6024533=3.7950934...
\]
and
\[
5-\bar t\bar p^{5/\bar t}=5-2^*(0.6024533)^{(5/2)}=4.43657195...
\]
Hence, as $\varepsilon\rightarrow 0$, we have
\[
\frac{z_H}{z^*}\rightarrow \frac{3.7950934}{4.43657195}\simeq 0.855411...
\]
\qed


\section*{Acknowledgements} The authors wish to acknowledge two anonymous reviewers for their helpful comments and remarks, which have significantly improved the rigor of the paper.

This material is based upon work supported by the National Science Foundation under Grant No. IIS-1909446.

\end{document}